# Improvement in Estimating Population Mean using Two Auxiliary Variables in Two-Phase Sampling

Rajesh Singh

Department of Statistics, Banaras Hindu University(U.P.), India

(rsinghstat@yahoo.com)

Pankaj Chauhan and Nirmala Sawan

School of Statistics, DAVV, Indore (M.P.), India

Florentin Smarandache

Department of Mathematics, University of New Mexico, Gallup, USA

(smarand@unm.edu)

**Abstract**

This study proposes improved chain-ratio type estimator for estimating population mean using some known values of population parameter(s) of the second auxiliary character. The proposed estimators have been compared with two-phase ratio estimator and some other chain type estimators. The performances of the proposed estimators have been supposed with a numerical illustration.

Key words: Auxiliary variables, chain ratio-type estimator, bias, mean squared error.

## 1. Introduction

The ratio method of estimation is generally used when the study variable Y is positively correlated with an auxiliary variable X whose population mean is known in advance. In the absence of the knowledge on the population mean of the auxiliary character we go for two-phase (double) sampling. The two-phase sampling happens to be a powerful and cost effective (economical) procedure for finding the reliable estimate in first phase sample for the unknown parameters of the auxiliary variable x and hence has eminent role to play in survey sampling, for instance, see Hidiroglou and Sarndal (1998).

Consider a finite population $U = (U_1, U_2, ......, U_N)$. Let y and x be the study and auxiliary variable, taking values $y_i$ and $x_i$ respectively for the $i^{th}$ unit $U_i$.

Allowing SRSWOR (Simple Random Sampling without Replacement) design in each phase, the two-phase sampling scheme is as follows:

(i) the first phase sample $s_{n'}$ ($s_{n'} \subset U$) of a fixed size $n'$ is drawn to measure only x in order to formulate a good estimate of a population mean $\bar{X}$,

(ii) Given $s_{n'}$, the second phase sample $s_n$ ($s_n \subset s_{n'}$) of a fixed size n is drawn to measure y only.

Let $\bar{x} = \frac{1}{n}\sum_{i \in s_n} x_i$, $\bar{y} = \frac{1}{n}\sum_{i \in s_n} y_i$ and $\bar{x}' = \frac{1}{n'}\sum_{i \in s_{n'}} x_i$.

The classical ratio estimator for $\bar{Y}$ is defined as

$$\bar{y}_r = \frac{\bar{y}}{\bar{x}}\bar{X} \qquad (1.1)$$

If $\bar{X}$ is not known, we estimate $\bar{Y}$ by two-phase ratio estimator

$$\bar{y}_{rd} = \frac{\bar{y}}{\bar{x}}\bar{x}' \qquad (1.2)$$

Some times even if $\bar{X}$ is not known, information on a cheaply ascertainable variable z, closely related to x but compared to x remotely related to y, is available on all units of the population. For instance, while estimating the total yield of wheat in a village, the yield and area under the crop are likely to be unknown, but the total area of each farm may be known from village records or may be obtained at a low cost. Then y, x and z are respectively yield, area under wheat and area under cultivation see Singh et.al.(2004).

Assuming that the population mean $\bar{Z}$ of the variable z is known, Chand (1975) proposed a chain type ratio estimator as

$$t_1 = \frac{\bar{y}}{\bar{x}}\left(\frac{\bar{x}'}{\bar{z}'}\right)\bar{Z} \tag{1.3}$$

Several authors have used prior value of certain population parameter(s) to find more precise estimates. Singh and Upadhyaya (1995) used coefficient of variation of z for defining modified chain type ratio estimator. In many situation the value of the auxiliary variable may be available for each unit in the population, for instance, see Das and Tripathi (1981). In such situations knowledge on $\bar{Z}, C_z$, $\beta_1(z)$ (coefficient of skewness), $\beta_2(z)$ (coefficient of kurtosis) and possibly on some other parameters may be utilized. Regarding the availability of information on $C_z$, $\beta_1(z)$ and $\beta_2(z)$, the researchers may be referred to Searls(1964), Sen(1978), Singh et.al.(1973), Searls and Intarapanich(1990) and Singh et.al.(2007). Using the known coefficient of variation $C_z$ and known coefficient of kurtosis $\beta_2(z)$ of the second auxiliary character z Upadhyaya and Singh (2001) proposed some estimators for $\bar{Y}$.

If the population mean and coefficient of variation of the second auxiliary character is known, the standard deviation $\sigma_z$ is automatically known and it is more meaningful to use the $\sigma_z$ in addition to $C_z$, see Srivastava and Jhajj (1980). Further, $C_z$, $\beta_1(z)$ and $\beta_2(z)$ are the unit free constants, their use in additive form is not much justified. Motivated with the above justifications and utilizing the known values of $\sigma_z$, $\beta_1(z)$ and $\beta_2(z)$, Singh (2001) suggested some modified estimators for $\overline{Y}$.

In this paper, under simple random sampling without replacement (SRSWOR), we have suggested improved chain ratio type estimator for estimating population mean using some known values of population parameter(s).

## 2. The suggested estimator

The work of authors discussed in section 1 can be summarized by using following estimator

$$t = \overline{y}\left(\frac{\overline{x}'}{\overline{x}}\right)\left(\frac{a\overline{Z}+b}{a\overline{z}'+b}\right) \qquad (2.1)$$

where $a(\neq 0)$, $b$ are either real numbers or the functions of the known parameters of the second auxiliary variable z such as standard deviation ($\sigma_z$), coefficient of variation ($C_z$), skewness ($\beta_1(z)$) and kurtosis ($\beta_2(z)$).

The following scheme presents some of the important known estimators of the population mean which can be obtained by suitable choice of constants a and b.

| Estimator | Values of | |
|---|---|---|
| | a | b |
| $t_1 = \overline{y}\left(\dfrac{\overline{x}'}{\overline{x}}\right)\left(\dfrac{\overline{Z}}{\overline{z}'}\right)$<br><br>Chand (1975) chain ratio type estimator | 1 | 0 |
| $t_2 = \overline{y}\left(\dfrac{\overline{x}'}{\overline{x}}\right)\left(\dfrac{\overline{Z}+C_z}{\overline{z}'+C_z}\right)$<br><br>Singh and Upadhyaya | 1 | $C_z$ |

| Estimator | | |
|---|---|---|
| (1995) estimator | | |
| $t_3 = \bar{y}\left(\dfrac{\bar{x}'}{\bar{x}}\right)\left(\dfrac{\beta_2(z)\bar{Z}+C_z}{\beta_2(z)\bar{z}'+C_z}\right)$  Upadhyaya and Singh (2001) estimator | $\beta_2(z)$ | $C_z$ |
| $t_4 = \bar{y}\left(\dfrac{\bar{x}'}{\bar{x}}\right)\left(\dfrac{C_z\bar{Z}+\beta_2(z)}{C_z\bar{z}'+\beta_2(z)}\right)$  Upadhyaya and Singh (2001) estimator | $C_z$ | $\beta_2(z)$ |
| $t_5 = \bar{y}\left(\dfrac{\bar{x}'}{\bar{x}}\right)\left(\dfrac{\bar{Z}+\sigma_z}{\bar{z}'+\sigma_z}\right)$  Singh (2001) estimator | 1 | $\sigma_z$ |
| $t_6 = \bar{y}\left(\dfrac{\bar{x}'}{\bar{x}}\right)\left(\dfrac{\beta_1(z)\bar{Z}+\sigma_z}{\beta_1(z)\bar{z}'+\sigma_z}\right)$  Singh (2001) estimator | $\beta_1(z)$ | $\sigma_z$ |
| $t_7 = \bar{y}\left(\dfrac{\bar{x}'}{\bar{x}}\right)\left(\dfrac{\beta_2(z)\bar{Z}+\sigma_z}{\beta_2(z)\bar{z}'+\sigma_z}\right)$ | $\beta_2(z)$ | $\sigma_z$ |

In addition to these estimators a large number of estimators can also be generated from the estimator t at (2.1) by putting suitable values of a and b.

Following Kadilar and Cingi (2006), we propose modified estimator combining $t_1$ and $t_i$ $(i=2,3,....,7)$ as follows

$$t_i^* = \alpha t_1 + (1-\alpha)t_i, \quad (i=2,3,....,7) \tag{2.2}$$

where $\alpha$ is a real constant to be determined such that MSE of $t_i^*$ is minimum and $t_i$ $(i=2,3,.....,7)$ are estimators listed above.

To obtain the bias and MSE of $t_i^*$, we write

$$\bar{y} = \bar{Y}(1+e_0), \quad \bar{x} = \bar{X}(1+e_1), \quad \bar{x}' = \bar{X}(1+e_1'), \quad \bar{z}' = \bar{Z}(1+e_2')$$

such that

$$E(e_0) = E(e_1) = E(e_1') = E(e_2') = 0$$

and

$$E(e_0^2) = f_1 C_y^2, \quad E(e_1^2) = f_1 C_x^2, \quad E(e_1'^2) = f_2 C_x^2$$

$$E(e_2'^2) = f_2 C_z^2, \quad E(e_0 e_1) = f_1 \rho_{xy} C_x C_y, \quad E(e_0 e_1') = f_2 \rho_{xy} C_x C_y$$

$$E(e_0 e_2') = f_2 \rho_{yz} C_y C_z, \quad E(e_1 e_1') = f_2 C_x^2, \quad E(e_1 e_2') = f_2 \rho_{xz} C_x C_z$$

$$E(e_1' e_2') = f_2 \rho_{xz} C_x C_z$$

where

$$f_1 = \left(\frac{1}{n} - \frac{1}{N}\right), \quad f_2 = \left(\frac{1}{n'} - \frac{1}{N}\right),$$

$$C_y^2 = \frac{S_y^2}{\bar{Y}^2}, \quad C_x^2 = \frac{S_x^2}{\bar{X}^2}, \quad C_z^2 = \frac{S_z^2}{\bar{Z}^2}$$

$$\rho_{xy} = \frac{S_{xy}}{S_x S_y}, \quad \rho_{xz} = \frac{S_{xz}}{S_x S_z}, \quad \rho_{yz} = \frac{S_{yz}}{S_y S_z}$$

$$S_y^2 = \frac{1}{(N-1)} \sum_{i \in U} (y_i - \bar{Y})^2, \quad S_x^2 = \frac{1}{(N-1)} \sum_{i \in U} (x_i - \bar{X})^2$$

$$S_z^2 = \frac{1}{(N-1)} \sum_{i \in U} (z_i - \bar{Z})^2, \quad S_{xy} = \frac{1}{(N-1)} \sum_{i \in U} (x_i - \bar{X})(y_i - \bar{Y})$$

$$S_{xz} = \frac{1}{(N-1)} \sum_{i \in U} (x_i - \bar{X})(z_i - \bar{Z}), \quad S_{yz} = \frac{1}{(N-1)} \sum_{i \in U} (y_i - \bar{Y})(z_i - \bar{Z}).$$

Expressing $t_i^*$ in terms of e's, we have

$$t_i^* = \bar{Y}(1+e_0)\left[\alpha(1+e_1')(1+e_1)^{-1}(1+e_2')^{-1} + (1-\alpha)(1+e_1')(1+e_1)^{-1}(1+\theta e_2')^{-1}\right] \quad (2.3)$$

where $\theta = \dfrac{a\bar{Z}}{a\bar{Z}+b}$ (2.4)

Expanding the right hand side of (2.3) and retaining terms up to second power of e's, we have

$$t_i^* \cong \overline{Y}[1 + e_0 - e_1 + e_1' - e_2'(\alpha + \theta - \alpha\theta)] \qquad (2.5)$$

or

$$t_i^* - \overline{Y} \cong \overline{Y}[e_0 - e_1 + e_1' - e_2'(\alpha + \theta - \alpha\theta)] \qquad (2.6)$$

Squaring both sides of (2.6) and then taking expectation, we get the MSE of the estimator $t_i^*$, up to the first order of approximation, as

$$MSE(t_i^*) = \overline{Y}^2[f_1 C_y^2 + f_3 C_x^2 + (\alpha + \theta - \alpha\theta)^2 f_2 C_z^2 - 2f_3 \rho C_y C_x - 2(\alpha + \theta - \alpha\theta) f_2 \rho C_y C_z]$$

$$(2.7)$$

where

$$f_3 = \left(\frac{1}{n} - \frac{1}{n'}\right).$$

Minimization of (2.7) with respect to $\alpha$ yield its optimum value as

$$\alpha_{opt} = \frac{K_{yz} - \theta}{1 - \theta} \qquad (2.8)$$

where

$$K_{yz} = \rho_{yz} \frac{C_y}{C_z}.$$

Substitution of (2.8) in (2.7) yields the minimum value of MSE ($t_i^*$) as –

$$min.MSE(t_i^*) = M_o = \overline{Y}^2[f_1 C_y^2 + f_3(C_x^2 - 2\rho_{yx} C_y C_x) - f_2 \rho_{yz}^2 C_y^2] \qquad (2.9)$$

3. Efficiency comparisons

In this section, the conditions for which the proposed estimator is better than $t_i (i = 1,2,....7)$ have been obtained. The MSE's of these estimators up to the order $o(n)^{-1}$ are derived as –

$$MSE(\overline{y}_{rd}) = \overline{Y}^2[f_1 C_y^2 + f_3(C_x^2 - 2\rho_{yx} C_y C_x)] \qquad (3.1)$$

$$MSE(t_1) = \overline{Y}^2[f_1 C_y^2 + f_2(C_z^2 - 2\rho_{yz} C_y C_z) + f_3(C_x^2 - 2\rho_{yx} C_y C_x)] \qquad (3.2)$$

$$\text{MSE}(t_2) = \overline{Y}^2 \left[ f_1 C_y^2 + f_2 (\theta_2^2 C_z^2 - 2\theta_2 \rho_{yz} C_y C_z) + f_3 (C_x^2 - 2\rho_{yx} C_y Cx) \right] \quad (3.3)$$

$$\text{MSE}(t_3) = \overline{Y}^2 \left[ f_1 C_y^2 + f_2 (\theta_3^2 C_z^2 - 2\theta_3 \rho_{yz} C_y C_z) + f_3 (C_x^2 - 2\rho_{yx} C_y Cx) \right] \quad (3.4)$$

$$\text{MSE}(t_4) = \overline{Y}^2 \left[ f_1 C_y^2 + f_2 (\theta_4^2 C_z^2 - 2\theta_4 \rho_{yz} C_y C_z) + f_3 (C_x^2 - 2\rho_{yx} C_y Cx) \right] \quad (3.5)$$

$$\text{MSE}(t_5) = \overline{Y}^2 \left[ f_1 C_y^2 + f_2 (\theta_5^2 C_z^2 - 2\theta_5 \rho_{yz} C_y C_z) + f_3 (C_x^2 - 2\rho_{yx} C_y Cx) \right] \quad (3.6)$$

$$\text{MSE}(t_6) = \overline{Y}^2 \left[ f_1 C_y^2 + f_2 (\theta_6^2 C_z^2 - 2\theta_6 \rho_{yz} C_y C_z) + f_3 (C_x^2 - 2\rho_{yx} C_y Cx) \right] \quad (3.7)$$

and

$$\text{MSE}(t_7) = \overline{Y}^2 \left[ f_1 C_y^2 + f_2 (\theta_7^2 C_z^2 - 2\theta_7 \rho_{yz} C_y C_z) + f_3 (C_x^2 - 2\rho_{yx} C_y Cx) \right] \quad (3.8)$$

where

$$\theta_2 = \frac{\overline{Z}}{\overline{Z} + C_z}, \quad \theta_3 = \frac{\beta_2(z)\overline{Z}}{\beta_2(z)\overline{Z} + C_z}, \quad \theta_4 = \frac{C_z \overline{Z}}{C_z \overline{Z} + \beta_2(z)}, \quad \theta_5 = \frac{\overline{Z}}{\overline{Z} + \sigma_z},$$

$$\theta_6 = \frac{\beta_1(z)\overline{Z}}{\beta_1(z)\overline{Z} + \sigma_z}, \quad \theta_7 = \frac{\beta_2(z)\overline{Z}}{\beta_2(z)\overline{Z} + \sigma_z}.$$

From (2.9) and (3.1), we have

$$\text{MSE}(\overline{y}_{rd}) - M_o = f_2 \rho_{yz}^2 C_y^2 \geq 0 \quad (3.9)$$

Also from (2.9) and (3.2)-(3.8), we have

$$\text{MSE}(t_i) - M_o = f_2 (\theta_i C_z - \rho_{yz} C_y)^2 \geq 0, \quad (i = 2, 3, \ldots, 7) \quad (3.10)$$

Thus it follows from (3.9) and (3.10) that the suggested estimator under optimum condition is always better than the estimator $t_i (i = 1, 2, \ldots, 7)$.

## 4. Empirical study

To illustrate the performance of various estimators of $\overline{Y}$, we consider the data used by Anderson (1958). The variates are

y : Head length of second son

x : Head length of first son

z : Head breadth of first son

$N = 25$, $\bar{Y} = 183.84 =$, $\bar{X} = 185.72$, $\bar{Z} = 151.12$, $\sigma_z = 7.224$, $C_y = 0.0546$, $C_x = 0.0526$, $C_z = 0.0488$, $\rho_{yx} = 0.7108$, $\rho_{yz} = 0.6932$, $\rho_{xz} = 0.7346$, $\beta_1(z) = 0.002$, $\beta_2(z) = 2.6519$.

Consider $n' = 10$ and $n = 7$.

We have computed the percent relative efficiency (PRE) of different estimators of $\bar{Y}$ with respect to usual estimator $\bar{y}$ and compiled in the table 4.1:

**Table 4.1: PRE of different estimators of $\bar{Y}$ with respect to $\bar{y}$**

| estimator | PRE |
| --- | --- |
| $\bar{y}$ | 100 |
| $\bar{y}_{rd}$ | 122.5393 |
| $t_1$ | 178.8189 |
| $t_2$ | 178.8405 |
| $t_3$ | 178.8277 |
| $t_4$ | 186.3912 |
| $t_5$ | 181.6025 |
| $t_6$ | 122.5473 |
| $t_7$ | 179.9636 |
| $t_i^*$ | 186.6515 |

## 5. Conclusion

We have suggested modified estimators $t_i^* (i = 2,3,....,7)$. From table 4.1, we conclude that the proposed estimators are better than usual two-phase ratio estimator $\bar{y}_{rd}$, Chand (1975) chain type ratio estimator $t_1$, estimator $t_2$ proposed by Singh and Upadhyaya (1995), estimators $t_i (i = 3,4)$ and than that of Singh (2001) estimators $t_i (i = 5,6,7)$. For practical purposes the choice of the estimator depends upon the availability of the population parameter(s).